\documentclass[10 pt]{amsart}
\def\a{\hbox{\bf a}} \def\b{\hbox{\bf b}} \def\c{\hbox{\bf c}}
  
\def\v{\hbox{\bf v}} \def\x{\hbox{\bf x}} \def\y{\hbox{\bf y}}
\def\w{\hbox{\bf w}}
 
\def\sa{\hbox{\bf{\tiny a}}} \def\sb{\hbox{\bf{\tiny b}}}
\def\sc{\hbox{\bf{\tiny c}}} 
 
 \def\sx{\hbox{\bf{\tiny x}}}
\def\sy{\hbox{\bf{\tiny y}}} \def\sw{\hbox{\bf{\tiny w}}}

\def\nn{{\mathbb N}} \def\zz{{\mathbb Z}} \def\rr{{\mathbb R}}

\def\bM{\bar{M}}

\def\no{\noindent} 
\def\di{\displaystyle} 
\def\bs{\bigskip\noindent}
\def\ms{\medskip\noindent}

\def\per{\hbox{\rm Per}} 
\def\sper{\hbox{\tiny Per}}

\def\tras{\hbox{\rm Trace}}

\def\BB{{\mathcal B}}  
\def\EE{{\mathcal E}} \def\LL{{\mathcal L}} \def\MM{{\mathcal M}}
  
\def\dLL{\partial{\mathcal L}}

\newtheorem{tmma}{\em Theorem}[section] 
\newtheorem{lemma}{\em Lemma}[section]

\newtheorem{corollary}{\em Corollary}[section] 
\newtheorem{remark}{\em Remark}[section]
 
\def\endpf{\hfill{$\Box$}}

\begin{document}

 
 
\title[Finite type Approximations]{
Finite type approximations of Gibbs measures\\ on sofic subshifts}
\author{J.-R.~Chazottes, L.~Ramirez and E.~Ugalde}

\begin{abstract} 
Consider a H{\"o}lder continuous potential $\phi$ defined on the full shift
$A^\nn$, where $A$ is a finite alphabet. Let $X\subset A^\nn$
be a specified sofic subshift. It is well-known that there is a unique
Gibbs measure $\mu_\phi$ on $X$ associated to $\phi$.
Besides, there is a natural nested sequence of subshifts of finite type $(X_m)$
converging to the sofic subshift $X$. To this sequence we can associate
a sequence of Gibbs measures $(\mu_{\phi}^m)$.
In this paper, we prove that these measures weakly converge at
exponential speed to $\mu_\phi$ ( in the classical
distance metrizing weak topology). We also establish a strong
mixing property (ensuring weak Bernoullicity) of $\mu_\phi$.
Finally, we prove that the measure-theoretic entropy of $\mu_\phi^m$
converges to the one of $\mu_\phi$ exponentially fast.
We indicate how to extend our results to more general subshifts
and potentials.
We stress that we use basic algebraic tools (contractive properties of
iterated matrices) and symbolic dynamics.

\end{abstract}

\maketitle
 

\bs\section{Introduction}\label{introduction}

\no Existence and uniqueness of equilibrium states/Gibbs measures associated to 
sufficiently regular potentials is established in the general context
of expansive homeomorphisms acting on a compact metric space satisfying
specification \cite{bowen75,hr92}.
This class of systems contains subshifts of finite types (coding Axiom
A diffeomorphisms) but more generally all specified subshifts like
topologically mixing sofic subshifts (on finite alphabets).

The usual way to prove existence and uniqueness is to
construct a sequence of elementary Gibbs measures (which are
atomic) and to argue that such a sequence must have an
accumulation point in the weak topology. Then one proves
that this accumulation point is unique.
In the particular case of subshifts of finite types and H{\"o}lder
continuous potentials, there is a complete theory of Gibbs measures
\cite{bowenbook}.

The point of view adopted here to study Gibbs measures on a specified
subshift $X$ is to approximate it by a nested sequence of subshifts of
finite type, $(X_m)$, in the sense of Hausdorff metric (there is a canonical way to do
this). This gives a natural sequences of Gibbs measures (finite-type
approximations) which converges weakly to a Gibbs measure whose properties we wish to analyze.

For the sake of definiteness, we assume that the given potential $\phi$ on
$A^\nn$ ($A$ is a finite alphabet) is H{\"o}lder continuous and $X\subset A^\nn$
is a specified sofic subshift. As we shall comment at the end of the paper,
we are not restricted to that situation. The two crucial properties on
which our method relies are specification and presence of magic words (see
definitions below). Sofic subshifts provide a natural class of subshifts with
such properties.

Our main result can be phrased as follows: The sequence of finite type approximations
$(\mu_{\phi}^m)$ defined on $(X_m)$ weakly converges, as $X_m\to X$, 
to a measure $\mu_\phi$ at an exponential speed. Then this measure
must be a Gibbs measure associated to $\phi$. Moreover, we prove a
strong mixing property (implying that $\mu_\phi$ is Bernoulli).
By a classical argument (Bowen), this implies uniqueness.
We also prove that the measure-theoretic
entropy $h(\mu_\phi^m)$ convergences to $h(\mu_\phi)$ exponentially
fast (as well as the relative entropy $h(\mu_\phi|\mu_\phi^m)$ to $0$).
We use and prove the fact that the topological pressure $P(\phi,X_m)$
converges to $P(\phi,X)$ exponentially fast.

We use two tools. The first one is algebraic
(contraction properties of iteration of primitive matrices
with respect to the projective metric); The second one is
symbolic dynamics.
All our constants have explicit expressions in terms of the `data' of
the problem, that is, the cardinality of the alphabet, the supremum
norm of the potential, its H{\"o}lder constants and the specification
length of the subshift $X$.

We would like to mention a related work to ours
due to Gurevich \cite{gurevich}. Therein the author deals with
measures of maximal entropy.
Informally speaking, he states some sufficient conditions on the way a subshift is
approximated in `entropy' by subshifts of finite type in order that the
corresponding sequence of measure of maximal entropy have a unique limit. The main
tool is graph theory.

\ms The paper is organized as follows.
In Section \ref{preliminary-notions} we record
basic definitions and notations.
Section \ref{main-results} contains our main results. 
Section \ref{technical-lemmas} is devoted to some
preparatory lemmas that we use for the proof of our
main results in Section \ref{proof-results}.
In Section \ref{final-remarks} we indicate some straightforward
generalizations of our results as well as examples.
We can indeed handle potential with polynomial variations (decaying
fast enough). Consequently, the exponential speeds mentioned above
become polynomial. We can also deal with more general specified subshifts (for instance,
non-sofic but specified $\beta$-shifts).

\bs\section{Preliminary notions}\label{preliminary-notions}


\bs \subsection{Symbolic dynamics}

Let $A$ be a finite alphabet.
For all integers $m,n$, $m\leq n$, in $\nn_0$ ($\nn_0:=\nn\cup \{0\}$),
we denote by $\a(m:n)$ the word $\a(m)\a(m+1)\cdots\a(n-1)\a(n)$ of
length $n-m+1$.
The distance
\begin{equation} 
d_A(\a,\b):=\exp\left(-\min\{n\geq 0:\a(-n:n)\neq\b(-n:n)\}
\right)
\end{equation}
makes the cartesian product $A^\nn$ a compact metric space.

\ms As usual, the {\it shift transformation} $T:A^\nn\to A^\nn$ is the map
such that $T\a(i)=\a(i+1)$.

\ms A {\it subshift} is a $T$--invariant compact set $X\subset
A^\nn$. The subshift $X$ is said to be of {\it finite type}, if
it is defined by a finite collection of admissible words, which
can be taken of the same length for the sake of simplicity (and
without loss of generality). So, the subshift of finite type
defined by the collection ${\LL}\subset A^{n+1}$ of {\it
admissible words}, is the compact set
\begin{equation}\label{definition-sft}
A_{\LL}:=\{\a\in A^\nn:\ \a(j:j+n)\in \LL \ \ \forall j\in\nn_0\}.
\end{equation}

\ms For a given subshift of finite type $X\subset A^\nn$, the {\it
order} of the subshift is the smallest integer $n\in\nn$ such that
$X$ is defined by a collection of admissible words of length $n$.

\ms A sequence $\a\in X$ is {\it periodic} of period $p\geq 1$ if
$T^p\a=\a$, and this is its minimal period if in addition
$T^k\a\neq \a$ whenever $0\leq k<p$. We will denote by $\per_p(X)$
the collection of all periodic sequences of period $p$ in $X$.

\ms For a general subshift $X\subset A^\nn$ and $n\geq 0$, the
collection of $X$--admissible words of length $m+1$ is the set
\begin{equation}\label{definition-admissible}
\LL_m(X):=\{\a(0:m):\ \a\in X\}.
\end{equation}

\ms {\it A sofic subshift} $X\subset A^\nn$ is a continuous
$T$--invariant image of a subshift of finite type. More precisely,
let $Y\subset A^\nn$ be a subshift of finite type, $B$ a finite
alphabet, and $\Pi:Y\to B^\nn$ a continuous map (with respect to
the distances $d_A$ and $d_B$), commuting with $T$. The image
$X=\Pi(Y)$, which in general is not of finite type, is a sofic
subshift.

\ms A more convenient way to characterize a sofic subshift is as
follows. Let $X\subset A^\nn$ be a subshift, and let
$\LL^*(X):=\cup_{n=0}^\infty\LL_n(X)$ be {\it the language}
defined by $X$.  For each $\a\in\LL^*(X)$ let $f(\a):=\left\{\b\in
\LL^*:\ \a\b\in\LL^*\right\}$ be the set of {\it followers} of
$\a$, and $p(\a):=\left\{\b\in \LL^*:\ \b\a\in\LL^*\right\}$ is
the set of {\it predecessors} of $\a$. The subshift $X$ is sofic
if $\{f(\a):\ \a\in\LL^*\}$ is a finite collection, in which case
$\{p(\a):\ \a\in\LL^*\}$ is finite as well~\cite{kitchensbook}.

\ms A word $\a\in\LL^{*}$ is a {\it magic word} for $X$ if $\b\in
p(\a)$ and $\c\in f(\a)$ implies $\b\a\c\in\LL^*$. It is a direct
consequence of the finiteness of the collection of followers that
every sofic subshift has a magic word (see~\cite[p.~148]{kitchensbook}).

\ms For a general subshift $X\subset A^\nn$ and a $X$--admissible
word $\a\in\LL_m(X)$, the set
\begin{equation}\label{definition-cylinder}
[\a]:=\{\b\in X:\ \b(0:m)=\a\}
\end{equation}
is the {\it cylinder} of length $m+1$ determined by $\a$. 

\ms The subshift $X$ is said to be {\it specified}, with
specification length $\ell\geq 1$ ($\ell=0$ means that we have a full shift),
if for each pair of $X$--admissible words $\a\in\LL_m(X)$ and $\b\in\LL_n(X)$, and
$k\geq \ell$, there exists a periodic sequence $\c$ of period
$m+n+k+2$, such that $\c(0:m)=\a$ and $(T^{m+k+1}\c)(0:n)=\b$.
Specification implies topological mixing and abundance of periodic
orbits in the sense that periodic orbits form a dense set in $X$.
See \cite{denkerbook} for more details on the
specification property.

\ms A notational remark: We shall use the symbols $\a,\b$, etc, both for
infinite sequences and finite words for convenience. To avoid any confusion
we shall always precise the nature of the $\a$'s or $\b$'s.


\bs\subsection{Gibbs measures}
The $\sigma$--field generated by the cylinders of $X\subset A^\nn$ 
coincides with the Borel
$\sigma$--field $\BB(X)$. The set $\MM(X)$ of Borel probability
measures in $\BB(X)$ is convex and compact in the weak topology.
The weak topology can be metrized with the distance
(see~\cite[p. 148]{waltersbook})
\begin{equation}\label{definition-weak}
D(\mu,\nu):= \sum_{m=0}^\infty 2^{-(m+1)}
\left(
\sum_{\sa\in\LL_m(X)} \left| \mu[\a]-\nu[\a]\right|
\right).
\end{equation}

\ms We denote by $\MM_T(X)$ the set of $T$--invariant probability measures
on $X$.

\ms A function $\phi:A^\nn\to \rr$ is H{\"o}lder continuous if for
some $\theta\in [0,1)$ and $C>0$, we have
$\max\{|\phi(\a)-\phi(\b)|: \a(0:m)=\b(0:m) \}|\leq C\theta^m$ for
all $m\geq 0$. The constant $\theta\in [0,1)$ is the {\it H{\"o}lder
exponent} of $\phi$. As usual, we shall call $\phi$ a potential.

\ms For $\phi:A^\nn\to\rr$ and $k\in \nn_0$ define
$S_k\phi:A^\nn\to\rr$ such that
\begin{equation}\label{definition-Smphi}
S_k\phi(\a)=\sum_{i=0}^k\phi\circ T^i(\a)\, .
\end{equation}

\ms Given a H{\"o}lder continuous potential $\phi$ and a subshift $X\subset A^\nn$,
$\mu\in\MM_T(X)$ is a {\it Gibbs measure} for the potential $\phi$
if there are constants $C=C(\phi,X) >0$ 
and $P(\phi,X)\in \rr$ such that
\begin{equation}\label{inequality-gibbs}
C^{-1} \leq \frac{\mu_\phi[\a(0:k)]}{\exp(S_k\phi(\a)-(k+1)P(\phi,X))}
\leq C
\end{equation}
for all $k\in \nn_0$.

\ms The constant $P(\phi,X)$ above, is the so called {\it topological
  pressure} of the potential $\phi$. For specified subshifts, 
it can be defined (see e.g. \cite{bowen75}) by the limit
\begin{equation}\label{defintion-pressure}
P(\phi,X):=\lim_{k\to\infty}\frac{1}{k+1}\log \left(
\sum_{\sa\in\sper_{k+1}(X)}\exp(S_k\phi(\a^*)) \right), 
\end{equation}
where $\a^*\in X$ is an arbitrary sequence in $[\a]$.

\ms For $X$ of finite type and $\phi$ H{\"o}lder continuous, there
exists a unique Gibbs measure $\mu_\phi$ (a proof of this fact can
be found in~\cite[p.~9~ff.]{bowenbook}, or \cite[ch.~5]{kellerbook}).
The existence and uniqueness of $\mu_\phi$ for general specified subshifts is a particular instance of
the {Theorem 2.5} in~\cite{hr92}.

\bs\section{Main results}\label{main-results}

\ms Let $X\subset A^\nn$ be a specified subshift.
The {\it finite type approximation of order $m$}, $m\in\nn$, 
to $X$ is the subshift of finite type
\begin{equation}\label{definition-ftapproximation}
X_m:=A_{\LL_m(X)}=\{\a\in A^\nn:\ \a(j:j+m)\in \LL_m(X),\; \forall
j\in \nn_0\},
\end{equation}
determined by the $X$--admissible words of length $m+1$. It is easy to
verify that the sequence of compact sets $\{X_m\}_{m\in \nn}$
converges in the Hausdorff metric to $X$ (you can find a definition
in~\cite[p. 111]{denkerbook}).

\ms On the finite type approximation $X_m$, the potential
$\phi:A^\nn\to\rr$ defines a unique Gibbs measure
$\mu_\phi^m\in\MM_T(X_m)$.
These measures will be used as {\it finite type approximations
of order $m$} of $\mu_\phi\in \MM_T(X)$.

\bs  For $m,p\in \nn$ let
$\EE_{(m,p)}\in\MM(X_m)$ be {\it the elementary Gibbs measure} with
support on $\per_{p+1}(X_m)$, such that
\begin{equation}\label{definition-elapproximation}
\EE_{(m,p)}[\b]:=
\frac{\exp(S_p\phi(\b))}{\sum_{\sa\in\sper_{p+1}(X_m)}
\exp(S_p\phi(\a))},
\end{equation}
for each $\b\in\per_{p+1}(X_m)$. 
We will use the fact \cite[p. 635]{KH} that each Gibbs measure 
$\mu_{\phi}^m$ can be obtained as a weak limit of the sequence of 
elementary Gibbs measures $\EE_{(m,p)}$, as $p\to\infty$.

\bs We have the following three main results, whose direct consequence
is the constructive proof of existence and uniqueness of Gibbs measures
on specified sofic subshifts, associated to H{\"o}lder continuous potentials. 

\begin{tmma}[\bf Speed of convergence of $\mu_\phi^m$]\label{theorem-elementary-convergence} 
Let $\phi:A^\nn\to\rr$ be a H{\"o}lder continuous potential, and
$X\subset A^\nn$ a {\bf\em  sofic specified subshift}. 
There exists an invariant measure $\mu^*\in\MM_T(X)$, a polynomial
$Q_{\rm FT}$ of degree $3$, and constants $\theta_{\rm FT}\in (0,1)$,
$m^*\in\nn$, 
satisfying 
\begin{equation}\label{speed}
D(\mu_\phi^m,\mu^*)\leq Q_{\rm FT}(m)\ \theta_{\rm FT}^{m},
\end{equation}
for all $m\geq m^*$.
\end{tmma}

\bs
\begin{tmma}[\bf `Gibbs property']\label{theorem-gibbs-inequality}
Under the hypotheses of Theorem~\ref{theorem-elementary-convergence}, there exists a constant
$C_g=C_g(X,\phi)>0$ such that
\begin{equation}
\label{gibbsproperty}
C_g^{-1}\leq \frac{\mu^*[\a(0:n)]}{\exp(S_{n}\phi(\a)-(n+1)P(\phi,
  X))}\leq C_g
\end{equation}
for each $n\in\nn_0$ and $\a\in X$. 
\end{tmma}

\bs
\begin{tmma}[\bf `Strong mixing']\label{theorem-mixing}
Under the hypotheses of Theorem~\ref{theorem-elementary-convergence}, 
there exists a polynomial $Q_\mu$ of degree 2, and $\theta_\mu\in (0,1)$ such that,
for all $\a,\b\in \LL^*(X)$ there exists $s^*:=s^*(\a,\b)$ satisfying
\begin{equation}\label{psi-mixing}
\left|\ \! \frac{\mu^*\left([\a]\cap T^{-s}[\b]\right)}{\mu^*[\a]\ \mu^*[\b]}- 1\ \!\!\right|\leq
Q_\mu(\sqrt{s})\ \theta_\mu^{\sqrt{s}} 
\end{equation}
for all $s\geq s^*$. 
\end{tmma}

\bs Combining the three previous theorems we get the following theorem.

\ms 
\begin{tmma}\label{theorem-principal}
Let $\phi:A^\nn\to\rr$ be a H{\"o}lder continuous potential, and
$X\subset A^\nn$ a {\bf\em  specified sofic subshift}. 
The weak limit  $\mu_\phi:=\lim_{m\to\infty}\mu^m_\phi$ is the unique 
Gibbs measure associated to the potential $\phi$, i.~e. the only 
$T$--invariant measure on $X$ satisfying \eqref{gibbsproperty}.
Moreover, the finite type approximations $\mu_\phi^m$ converge
exponentially fast to $\mu_\phi$ in the sense of \eqref{speed}
and $\mu_\phi$ is mixing in the sense of \eqref{psi-mixing} and Bernoulli.
\end{tmma}

\ms 
\begin{proof}
Theorems~\ref{theorem-elementary-convergence}
and~\ref{theorem-gibbs-inequality}
ensure the existence of a measure satisfying the inequalities \eqref{gibbsproperty}
and having exponentially fast converging finite type approximations.
To prove uniqueness, we can follow the last part of the proof of
Theorem 1.16 in~\cite{bowenbook}.
The mixing property \eqref{psi-mixing} implies weak Bernoullicity,
see e.g. \cite[p. 169]{shieldsbook}. The theorem is proved.

\end{proof}

\begin{remark}
All constants appearing in the above theorems, including the
coefficients of the polynomials, have explicit (but somewhat tedious)
expressions in terms of the data of the problem, that is,
$\# A$, $\parallel\!\!\phi\!\!\parallel$, $C$, $\theta$ (H{\"o}lder
condition) and $\ell$ (the specification length).
These expressions are given in the proofs.
\end{remark}

We end this section with the following theorem on speed of convergence
of the entropy $h(\mu_\phi^m)$ to $h(\mu_\phi)$, and the
relative entropy $h(\mu_\phi | \mu_\phi^m)$ to $0$.

\begin{tmma}\label{entropy-approximation}
Under the hypotheses of Theorem~\ref{theorem-elementary-convergence}, there exist constants
$C_h>0$, $C_P>0$, $0<\theta_h<1$ and $0<\theta_P<1$, such that
\begin{equation}\label{entropy}
|h(\mu_\phi)-h(\mu_\phi^m)|\leq C_h\ \theta_h^m
\end{equation}
\begin{equation}\label{relative-entropy}
h(\mu_\phi | \mu_\phi^m)\leq \frac{C_P}{1-\theta_P}\ \theta_P^m\,.
\end{equation}
\end{tmma}

We refer the reader to \cite{waltersbook} for details on
entropy of invariant measures. The appendix at the end of the
paper contains the necessary informations on entropy and
relative entropy regarding our context.

To the best of our knowledge, Theorems
\ref{theorem-elementary-convergence}-\ref{theorem-gibbs-inequality}-\ref{theorem-mixing}
and \ref{entropy-approximation} are new. The first three ones imply
existence and uniqueness of $\mu_\phi$. The only known mixing property 
for this measure is the usual mixing property (which does not assure
Bernoullicity). This mixing much less stronger than \eqref{psi-mixing}
which implies Bernoullicity.

\bs\section{Technical lemmas}\label{technical-lemmas}

\ms In this section we establish some technical lemmas
needed to prove theorems of Section \ref{main-results}.
We shall use some results coming from the theory of
primitive matrices, as well as some elementary facts about weak
distance between measures. The Appendix contains these results and
some related notions. From now on we assume known those results and
notions, as well as the notations established there.

\bs {\bf Notations}. From now on, an expression of the type $a=c^{\pm 1}$
stands for the inequalities $c^{-1}\leq a\leq c$. Similarly
$a={\pm}c$ stands for $-c\leq a\leq c$. By extension, $a=\exp(\pm b)$
will stand for $\exp(-b)\leq a \leq \exp(b)$. 

\bs Given a H{\"o}lder continuous potential $\phi: A^\nn\to \rr$, for
each $n\in\nn_0$ we define the finite range potential $\phi^n:A^{n+1}\to\rr$ such that
\begin{equation}\label{definition-phiapproximation}
\phi^n(\a)=\max\{\phi(\b):\ \b\in[\a]\}.
\end{equation}

\bs For $n\geq m$ let $\LL_{(m,n)}:=\LL_n(X_m)$ be the set of
$X_m$--admissible words of length $n+1$, which of course contains
$\LL_n:=\LL_n(X)$. Let us define the transfer matrix
$M_{(m,n)}:\LL_{(m,n)}\times\LL_{(m,n)}\rightarrow\rr^+$ such that
\begin{equation} \label{definition-transfermatrix}
M_{(m,n)}(\a,\b)= \left\{
\begin{array}{ll}
\exp(\phi^{n+1}(\a\b(n)) & \text{if }\ \a(1:n)=\b(0:n-1), \\ 0 &
\text{otherwise.}
\end{array}
\right.
\end{equation}

\ms For a specified subshift $X$, the matrix $M_{(m,n)}$ is primitive
with primitivity index $\ell+n+1$, and
has a unique maximal eigenvalue $\rho_{(m,n)}:=\rho(M_{(m,n)})$.
Associated to $\rho_{(m,n)}$ there are unique normalized right and left eigenvectors
$\v_{(m,n)}:=\v_{M_{(m,n)}}$ and $\w_{(m,n)}:=\w_{M_{(m,n)}}$.

\ms The elementary measure $\EE_{(m,p)}$ can be expressed in 
term of the transfer matrices $M_{(m,n)}$ as follows.

\ms For $p > m \geq n$, and $\a\in \LL_{(m,n)}$, we have
\begin{equation}\label{definition-elementarymatrix}
\EE_{(m,p)}[\a]=
\frac{M_{(m,n)}^{p+1}(\a,\a)}{\tras\left(M_{(m,n)}^{p+1}\right)}
\times \exp(\pm 2(p+1)C\theta^{n+1}).
\end{equation}

\ms Now, given $n > m$, for each each $\a\in\LL_{(m,n)}$ 
define $R^{\sa}, L_{\sa}: \LL_{(m,n)}\to\rr^+ $ be such that
\begin{equation}\label{definition-brackets}
R^{\sa}(\b)=M_{(m,n)}^{\ell+n+1}(\b,\a) \ \text{ and }
L_{\sa}(\b)=M_{(m,n)}^{\ell+n+1}(\a,\b).
\end{equation}
Note that these vectors are positive.

\ms We are able to give a uniform estimate of the values of elementary
measures on cylinders, by using Corollary~\ref{corollary-powermatrix}.

\bs
\begin{lemma}\label{lemma-elementarymatrix}
Let $X\subset A^\nn$ be a specified subshift with specification length
$\ell\geq 0$, $\phi:A^\nn\to\rr$ a H{\"o}lder continuous potential with
constant $C>0$ and exponent $\theta\in(0,1)$.  There are constants $C_\EE > 0$
and $\theta_\EE \in (0,1)$ such that, for all integers $m,n,p$,
such that $m\leq n$, $(n+1)(n+\ell+1)\leq p$ and $\a\in\LL_{(m,n)}$, we have
\[
\EE_{(m,p)}[\a]=\w_{(m,n)}(\a)\v_{(m,n)}(\a)\times
\exp(\pm (p+1)\ C_\EE {\theta_\EE}^{n+1}).
\]
\end{lemma}

\ms 
\begin{proof} 
For each $m\leq n$ let $\tau_{(m,n)}$ be the Birkhoff contraction
coefficient of $M_{(m,n)}^{n+\ell+1}$. Let $\bM:=M_{(m,n)}^{n+\ell+1}$
and $\Gamma:=\Gamma(\bM)\in (0,1]$, where $\Gamma$ is defined in
\eqref{defintion-gamma} in the Appendix.
According to Theorem~\ref{theorem-seneta}, we have $\left(1-\tau_{(m,n)}\right)^{-1}
=(1+\Gamma)/(2\Gamma) \leq \Gamma^{-1}$.

\ms On the other hand we have
\begin{eqnarray*}
\Gamma^{-1}&=& \min_{\sa,\sb,\sa',\sb'\in\LL_{(m,n)}}
        \left(\frac{\bM(\a,\b)\bM(\a',\b')}{\bM(\a,\b')\bM(\a',\b))}\right)^{-1/2}\\
                                &=&
                                \max_{\sa,\sb,\sa',\sb'\in\LL_{(m,n)}}
     \left(\frac{\bM(\a,\b)\bM(\a',\b')}{\bM(\a,\b')\bM(\a',\b))}\right)^{1/2}.
\end{eqnarray*}
Now, for arbitrary $\a,\b,\a',\b'\in\LL_{(m,n)}$ we have
$$
\frac{\bM(\a,\b)\bM(\a',\b')}{\bM(\a,\b')\bM(\a',\b)} \leq
$$
$$
\frac{\di
\sum_{\sc=\sa\sx\sb\in\LL_{(m,\ell+2n+1)}}\!\!\!\!\!\!\!\!\!\!\!\!\!
\exp(S_{n+\ell}\phi^{n+1}(\c))
}{\di
\min_{\sc=\sa\sx\sb\in\LL_{(m,\ell+2n+1)}}\!\!\!\!\!\!\!\!\!\!\!\!\!
\exp(S_{n+\ell}\phi^{n+1}(\c))}
\times
\frac{\di\sum_{\sc'=\sa'\sx\sb'\in\LL_{(m,\ell+2n+1)}}\!\!\!\!\!\!\!\!\!\!\!\!\!
\exp(S_{n+\ell}\phi^{n+1}(\c'))
}{
\di\min_{\sc'=\sa'\sx\sb'\in\LL_{(m,\ell+2n+1)}}\!\!\!\!\!\!\!\!\!\!\!\!\!
\exp(S_{n+\ell}\phi^{n+1}(\c'))}\\ 
$$
$$                        
\leq \left(\left( \# A\ e^C\right)^\ell\times
e^{\Lambda \theta}\right)^2\; .
$$
Hence 
\[
\frac{1}{1-\tau_{(m,n)})}\leq K_0:=
\left(e^C \# A\right)^\ell\times e^{\Lambda\theta}.
\]
i.e. $\tau_{(m,n)} \leq 1-K_0^{-1} < 1$.

\bs For each $m\leq n$ let $d_{(m,n)}$ be the projective distance on the
simplex of dimension $\#\LL_{(m,n)}$, and $F_{(m,n)}$ the
transformation defined on the simplex by the transition matrix
$M_{(m,n)}$. Note that
\[ 
d_{(m,n)}\left(M_{(m,n)} R^{\sa},R^{\sa}\right)= \log\left(
\frac{\max_{\sb}M_{(m,n)}^{n+\ell+2}(\b,\a)/M_{(m,n)}^{n+\ell+1}(\b,\a)
}{ \min_{\sb}M_{(m,n)}^{n+\ell+2}(\b,\a)/M_{(m,n)}^{n+\ell+1}(\b,\a)}
\right).
\] 
We have
\begin{eqnarray*}
\frac{M_{(m,n)}^{n+\ell+2}(\b,\a)}{M_{(m,n)}^{n+\ell+1}(\b,\a)} &=&
        \frac{\di\sum_{\sc=\sb\sx\sa\in\LL_{(m,2n+\ell+2)}}\exp(S_{n+\ell+1}\phi^{n+1}(\c))
        }{\di\sum_{\sc'=\sb\sy\sa\in\LL_{(m,2n+\ell+1)}}\exp(S_{n+\ell}\phi^{n+1}(\c'))}\\
                &=& 
        e^{\pm||\phi||}\left(\# A \ e^C\right)^{\pm(\ell+1)}e^{\Lambda\theta}\,.
\end{eqnarray*}
where $\parallel\!\!\phi\!\!\parallel:=\max\{|\phi(\a)|:\ \a\in A^{\nn}\}$.
From this we get
\[
\max_{\a\in\LL_{(m,n)}}d_{(m,n)}\left(F_{(m,n)}
\left(R^{\sa}\right),R^{\sa}\right) \leq K_1
\]
with 
\[
K_1:=2 \left((\ell+1)(\log(\# A)+C)+\Lambda\theta +||\phi|| \right)\, .
\]

\bs Finally, with (\ref{definition-brackets}), 
inequalities~(\ref{definition-elementarymatrix}) may be rewritten as
\[
\EE_{(m,p)}[\a]= \frac{L_{\sa}^\dag M_{(m,n)}^{p+1-(n+\ell+1)}R^{\sa}
}{\di \sum_{\sb\in\LL_{(m,n)}}L_{\sb}^\dag M_{(m,n)}^{p+1-(n+\ell+1)}
R^{\sb}} \times \exp(\pm 2(p+1)C\theta^{n+1})\, .
\]
Then, using Corollary~\ref{corollary-powermatrix}, we have
$$
\EE_{(m,p)}[\a] = \frac{\w_{(m,n)}^\dag R^{\sa}\ 
L_{\sa}^\dag\v_{(m,n)}}{\di \sum_{\sb\in\LL_{(m,n)}} \w_{(m,n)}^\dag \
R^{\sb} L_{\sb}^\dag\v_{(m,n)}} \quad \times 
$$
$$
\exp\left(\pm 
2\left( (p+1)C\theta^{n+1} + 
        K_0K_1(n+\ell+1)\left(1-K_0^{-1}\right)^{\lfloor \frac{p+1}{n+\ell+1} \rfloor} 
 \right)\right)
$$
where $K_0$ and $K_1$ are given above.
On the other hand we have 
\begin{eqnarray*}
L_{\sa}^\dag \v_{(m,n)}&=&
    \left(M_{(m,n)}^{n+\ell+1}\v_{(m,n)}\right)(\a)=
    \rho_{(m,n)}^{n+\ell+1}\v_{(m,n)}(\a),\\ 
\w_{(m,n)}^{\dag}R^{\sa}&=&
    \left(\w_{(m,n)}^{\dag}M_{(m,n)}^{n+\ell+1}\right)(\a)=
    \rho_{(m,n)}^{n+\ell+1}\w_{(m,n)}(\a),
\end{eqnarray*}
and $\w_{(m,n)}^{\dag}\v_{(m,n)}=1$. Then, taking into account 
that $p+1\geq (n+1)(n+\ell+1)$
and $m\leq n$, we obtain
\[
\EE_{(m,p)}[\a] = \w_{(m,n)}(\a)\v_{(m,n)}(\a)\times 
\exp\left(\pm (p+1)\ C_\EE\ {\theta_\EE}^{n+1}\right)
\]
with $C_\EE:=2(C + K_0 K_1)$ and $\theta_\EE:=\max\left(1-K_0^{-1}, \theta\right)$.
The lemma is proved.
\end{proof}

\bs
\begin{lemma}\label{lemma-elementaryelementary}
Let $X\subset A^\nn$ be a specified sofic subshift, with specification
length $\ell\geq 1$, and $\phi:A^{\nn}\to\rr$ a H{\"o}lder continuous potential
with constant $C>0$ and exponent $\theta\in (0,1)$. Then there are
constants $m_X\in\nn$, $C_X >0$ and $\theta_X\in(0,1)$, such that for
$m_X \leq  m \leq p$  
\[
1-(p+1)C_X\theta_X^m \leq \frac{\EE_{(m,p)}[\a]}{\EE_{(m+1,p)}[\a]}
\leq 1
\]
for each $\a\in \per_{p+1}(X_{m+1})$.
\end{lemma}

\bs \begin{proof} First note that
$$
\frac{\EE_{(m,p)}[\a]}{\EE_{(m+1,p)}[\a]} =  \frac{
\sum_{\sb\in \sper_{p+1}(X_{m+1})}e^{S_p\phi(\sb)} }{
\sum_{\sb\in\sper_{p+1}(X_{m}) } e^{S_p\phi(\sb)} } = 
$$ 
$$
1-\frac{\sum_{\sb\in\sper_{p+1}(X_m\setminus X_{m+1})}
e^{S_p\phi(\sb)}}{\sum_{\sb\in\sper_{p+1}(X_{m})}e^{S_p\phi(\sb)}}\geq 
$$
$$
1-(p+1)\, \frac{\sum_{\sb\in\sper_{p+1}(\partial X_m)}
\exp(S_p\phi(\b))
}{\sum_{\sb\in\sper_{p+1}(X_{m})}\exp(S_p\phi(\b))}
$$
where $\partial X_m:=\{\a\in X_m:\ \a(0:m+1)\not\in \LL_{m+1}\}$.

\ms Let $\dLL_{m}:=\LL_{(m,m+1)}\setminus\LL_{m+1}$.
Using specification property we obtain
\[
\frac{\sum_{\sb\in\sper_{p+1}(\partial X_m)} e^{S_p\phi(\sb)}
}{\sum_{\sb\in\sper_{p+1}(X_{m})}e^{S_p\phi(\sb)}}\leq
\frac{\sum_{\sa\in\dLL_{m}} e^{S_{m+1}\phi(\sa^*)}
}{\sum_{\sa\in\LL{(m,m+1)}}e^{S_{m+1}\phi(\sa^*)}}\times \left(\#A
e^{2\parallel\phi\parallel}\right)^{\ell}e^{4\Lambda}
\]
for any $\a^*\in[\a]$. We will prove that the quotient
$$
{\left(\sum_{\sa\in\dLL_{m}}
\exp(S_{m+1}\phi(\a^*))\right)\Big/
\left(\sum_{\sa\in\LL{(m,m+1)}}\exp(S_{m+1}\phi(\a^*))\right)}
$$
is exponentially small with $m$. 
This is the point at which we use the existence of
magic words.

\ms Fix a magic word $\w\in\LL_k$ with $k\geq \ell+1$. This is
always possible since for a magic word $\a\in\LL^*$, the
concatenated word $\a \b$ is again magic, for any $\b\in f_X(\a)$
($f_X(\a)$ is the set of followers of $\a$, which contains arbitrary
long words).
Let $m\geq 2k(k+\ell)$, so
that $\lfloor (m+1)/(k+\ell+1)\rfloor \geq m/k$ (we will use this
condition at the final step of the proof). Note that if
$\a\in\dLL_m$, then $\a(i:i+k)\neq \w$ for each $1\leq i\leq m-k$.
This is because if $\a(i:i+k)=\w$ then $\a(0:i+k),\a(i:m+1)\in
\LL^*(X)$, implying that $\a\in\LL^*(X)$ which contradicts the
hypothesis.

\ms Letting $q:=k+\ell+1$ define
\[
\dLL_m^{\sw}:=\{\a\in \dLL_m:\ \a(jq:jq+k)\neq \w,\ 0\leq j\leq
\lfloor (m+1)/q\rfloor-1 \}.
\]
It is clear that $\dLL_m\subset \dLL_m^{\sw}$. Define also
\[
\epsilon^{\sw}:=\frac{
\exp(S_k\phi(\b^-))}{ \sum_{\sb\in\LL_k\setminus\{\sw\}
}\exp(S_k\phi(\b^+))}\times \left(\#A\ e^{2||\phi||}\right)^{-\ell},
\]
where for each $\b\in\LL_k$, the sequences $\b^-,\b^+\in[\b]$
are such that $S_k\phi(\b^-)=\min_{\sb^*\in[\sb]}S_k\phi(\b^*)$
and $S_k\phi(\b^+)=\max_{\sb^*\in[\sb]}S_k\phi(\b^*)$.

\ms Let $r:=\lfloor (m+1)/q\rfloor$. For each
$\omega\in\{0,1\}^r$ define
\[
\LL_{(m,m+1)}^{\omega}:=\left\{\a\in\LL_{(m,m+1)}:\ \a(jq:jq+k)=\w
\ \text{if and only if}\ \omega(j)=1\right\}.
\]
It is clear that the collection $\left\{\LL_{(m,m+1)}^{\omega}:\
\omega\in\{0,1\}^r\right\}$ is a partition of $\LL_{(m,m+1)}$.
Now, it follows from the specification property that for each
$\omega\in\{0,1\}^r$
\[
\sum_{\sb\in\LL_{(m,m+1)}^\omega}\exp(S_{m+1}\phi(\b^*))\geq
\left(\epsilon^{\sw}\right)^{|\omega|_1}\times
\sum_{\sb\in\dLL_{m}^{\sw}}\exp(S_{m+1}\phi(\b^-))
\]
where, as before, $\b^-\in[\b]$ minimizes $S_{m+1}\phi$, and
$|\omega|_1:=\sum_{j=0}^{r-1}\omega(i)$. From the previous
inequality we readily derive
\[
\sum_{\sb\in\LL_{(m,m+1)}}\exp(S_{m+1}\phi(\b^*))\geq
\left(1+\epsilon^{\sw}\right)^r\times
\sum_{\sb\in\dLL_{m}^{\sw}}\exp(S_{m+1}\phi(\b^{-}))\,.
\]

\ms Finally,
\[
\frac{\sum_{\sa\in\dLL_{m}} \exp(S_{m+1}\phi(\a^*))
}{\sum_{\sa\in\LL{(m,m+1)}}\exp(S_{m+1}\phi(\a^*))}\leq
\frac{\sum_{\sa\in\dLL_{m}^{\sw}} \exp(S_{m+1}\phi(\a^*))
}{\sum_{\sa\in\LL{(m,m+1)}}\exp(S_{m+1}\phi(\a^*))}\\
\leq
\left(1+\epsilon^{\sw}\right)^{-r}\,.
\]
Since $m\geq 2k(k+\ell)$ then
$\left(1+\epsilon^{\sw}\right)^{-r}\leq
\left(1+\epsilon^{\sw}\right)^{-m/k}$, and the result follows with
$$
C_X:=\left(\#A e^{2\parallel\phi\parallel}\right)^\ell e^{4\Lambda},\;\;
\theta_X:=\left(1+\epsilon^{\sw}\right)^{-1/k},\;\;m_X=2k(k+\ell)\, .
$$
The lemma is proved.
\end{proof}

The following lemma has its own interest.

\bs
\begin{lemma}\label{lemma-pressure}
Let $X\subset A^\nn$ be a specified sofic subshift, with specification
length $\ell$. Let $\phi:A^{\nn}\to\rr$ be a H{\"o}lder continuous potential
with constant $C>0$ and exponent $\theta\in (0,1)$. Then there are constants
$m_P\in\nn$, $C_P >0$ and $\theta_P\in(0,1)$, such that 
\[
0 \leq P(\phi,X_m)-P(\phi,X_{m+1})\leq C_P\theta_P^m
\]
for all $m\geq m_P$.
\end{lemma}

\begin{proof} Proceeding as in the proof of the previous lemma, we obtain
\begin{eqnarray*}
0&\leq& \frac{1}{p+1}\log\left(
\frac{ \sum_{\sa\in \sper_{p+1}(X_{m})}\exp(S_{p+1}\phi(\a)) }
{ \sum_{\sa\in\sper_{p+1}(X_{m+1})} \exp(S_{p+1}\phi(\a)) }\right) \\
&\leq& \frac{1}{p+1}\log\left(1+\frac{(p+1)C_X\theta_X^m}{1-(p+1)C_X\theta_X^m}\right)\\
&\leq& \frac{C_X\theta_X^m}{1-(p+1)C_X\theta_X^m}
\end{eqnarray*} 
for $m\geq m_X$.

\bs To make use of the previous inequality, we need to know the speed of convergence of
\[
\frac{1}{p+1}\log\left(
\sum_{\sa\in\sper_{p+1}(X_m)}\exp(S_p\phi(\a))\right)\quad\quad\textup{to}\quad\quad
P(\phi,X_m).
\]

\ms Some computations like the ones done to prove
Lemma~\ref{lemma-elementarymatrix} give 
we obtain 
$$
\sum_{\sa\in\sper_p(X_m)}\exp(S_{p+1}\phi(\a))
=\tras\left(M_{(m,n)}^{p+1}\right)\times\exp(\pm (p+1)C\theta^{n+1})=
$$
$$
\left(\sum_{\sb\in\LL_{(m,n)}} \w_{(m,n)}^\dag R^{\sb} L_{\sb}^\dag\v_{(m,n)}\right) 
\times\rho_{(m,n)}^{p+1-2(n+\ell+1)}\times
\exp\left(\pm C_{\EE}(p+1) \theta_{\EE}^{n+1}\right)=
$$
$$
\rho_{(m,n)}^{p+1}\times\exp\left(\pm C_{\EE}(p+1) \theta_{\EE}^{n+1}\right)            
$$
for each $m < n$ and $(n+1)(n+\ell+1)\leq p$. Therefore
\[
\frac{1}{p+1}\log\left(\sum_{\sa\in\sper_{p+1}(X_m)}\exp(S_p\phi(\a))\right)=
\log\rho_{(m,n)}\pm C_{\EE}\theta_{\EE}^{n+1}.
\]

\ms Let us now prove that $\left\{\rho_{(m,n)}\right\}_{n > m}$
converges exponentially fast.
By definition, the limit has to be equal to $\exp(P(\phi,X_m))$.

\ms Let us define $N:\LL_{(m,n+1)}\times \LL_{(m,n+1)} \to \rr^+$ such that 
\[
N(\a,\b)=\left\{\begin{array}{ll}
\exp(\phi^{n+1}(\a)) & \text{if }\ \a(1:n+1)=\b(0:n)\\
                0  & \text{otherwise.}\end{array}\right.
\] 
Note that $M_{(m,n)}=N\exp(\pm C\theta^{n+1})$ coordinate-wise 
and $\rho_{(m,n)}=\exp(\pm C\theta^{n+1})\rho_N$. 
This can be easily derived from Corollary~\ref{corollary-powermatrix},
taking into account that 
$\rho_M=\lim_{n\to\infty}\left(\y^\dag M^n\x\right)^{1/n}$ for a primitive matrix $M$,
and arbitrary positive vectors $\x,\y$. 
Let $\v:\LL_{(m,n)}\to\rr^+$ such that 
$\v(\a)=\exp(\phi^{n+1}(\a))\times\v_{(m,n)}(\a(1:n+1))$, we have 
\begin{eqnarray*}
(N\v)(\a) &=& \exp(\phi^{n+1}(\a)) (M_{(m,n)}\v_{(m,n)})(\a(1:n+1))\\
          &=& \exp(\phi^{n+1}(\a))\rho_{(m,n)}\v_{(m,n)}(\a(1:n+1))
           = \rho_{(m,n)}\x(\a).
\end{eqnarray*} Hence, $\v$ is a positive eigenvector for the
matrix $N$, associated to the positive eigenvalue
$\rho_{(m,n)}$. Since $N$ primitive, Corollary~\ref{corollary-powermatrix} implies that
$\rho_N=\rho_{(m,n)}$, therefore $\rho_{(m,n+1)}=\exp(\pm C \theta^{n+1})\rho_{(m,n)}$.
>From this we obtain,
\[\frac{\rho_{(m,n)}}{\exp(P(\phi,X_m))}=
\exp\left(\pm \Lambda\theta^{n+1}\right)\, .
\]

\ms Since $X_m\supset X_{m+1}$, then $P(\phi,X_m)\geq
P(\phi,X_{m+1})$.
The previous computations imply on the other hand that
\[
P(\phi,X_m)-P(\phi,X_{m+1})\leq
\]
\[
\frac{C_X\theta_X^m}{1-((m+2)(m+\ell+2)+1)C_X\theta_X^m} +
C_{\EE}\theta_{\EE}^{m+2} + \Lambda\theta^{m+2},
\]
for $m\geq m_X$, by taking $n=m+1$ and $p=(n+1)(n+\ell+1)$. 
Thus, the lemma follows with 
\begin{eqnarray*}
\theta_P&:=&\max(\theta,\theta_X,\theta_\EE),\\ 
     C_P&:=&2C_X+C_\EE\theta_\EE^2+\Lambda\theta^2
\end{eqnarray*}
and $m_P:=\max(m_X, m_0)$, with $m_0$ such that
$2C_X((m+2)(m+\ell+2)+1)\theta_X^m\leq 1$ for all $m\geq m_0$. 

\end{proof}

\bs\section{Proof of the main results}\label{proof-results}

\ms This section is devoted to the proof of Theorems \ref{theorem-elementary-convergence}, \ref{theorem-gibbs-inequality},
\ref{theorem-mixing} and \ref{entropy-approximation}.  

\bs\subsection{Proof of Theorem~\ref{theorem-elementary-convergence}}

\ms Lemma~\ref{lemma-elementarymatrix} implies that
\[
\EE_{(m,p+1)}[\a]=\exp(\pm 2(p+2)C_\EE\theta^{\sqrt{p}-(\ell/2+1)})\ \EE_{(m,p)}[\a],
\] 
for each $\a\in \cup_{k=1}^{\lfloor \sqrt{p}-(\ell/2 +1)\rfloor}\LL_{(m,k)}$. Then 
Lemma~\ref{lemma-marginalsratio} applies, and we obtain  
\begin{eqnarray*}
D(\EE_{(m,p)},\EE_{(m,p+1)}) &\leq& 
\left(\exp\left( 2C_{\EE}(p+2)\ \theta_{\EE}^{\sqrt{p}-(\ell/2+1)}\right)-1\right)
                                                      +2^{(\ell/2+1)-\sqrt{p}}\\
                              &\leq&
 4C_{\EE}(p+2)\ \theta_{\EE}^{\sqrt{p}-(\ell/2+1)} + 2^{(\ell/2+1)-\sqrt{p}}
\end{eqnarray*}
for each $p\geq \max(p_0,(m+2)(m+\ell+2))$, with 
\[p_0:=\min\left\{p\in\nn:\ 2C_{\EE}(k+2)\theta_{\EE}^{\sqrt{k}-(\ell/2+1)}\leq 1
\ \text{for all}  \ k\geq p\right\}.
\]
Since $\sum_{p=0}^{\infty}(p+2)\ \theta_{\EE}^{\sqrt{p}}<\infty$, there exists
a  limit measure $\mu^m:=\lim_{p\to\infty}\EE_{(m,p)}$ belonging to $\MM_T(X_m)$. 
The convergence is such that
\[
D(\mu^m,\EE_{(m,p)})\leq 4C_\EE\theta_\EE^{-(\ell/2+1)}
Q(\sqrt{p})\theta_\EE^{\sqrt{p}} + 2^{(\ell/2+3)}2^{-\sqrt{p}}(\sqrt{p}+3)
\]
for each $p\geq \max(p_0,(m+2)(m+\ell+2))$. Here
\[
Q(x):=-\frac{2x(x^2+3)}{\log(\theta_\EE)}+\frac{6(x^2+1)}{\log^2(\theta_\EE)}
-\frac{12x}{\log^3(\theta_\EE)}+\frac{12}{\log^4(\theta_\EE)}.
\]

\ms Let us now prove that the limiting measure $\mu^m$ coincides
with the unique Gibbs measure $\mu_\phi\in\MM_T(X_m)$.
From the specification property we can derive the inequalities
\begin{eqnarray*}
\EE_{(m,p)}[\a]&=&
\exp\left(S_n\phi(\a^*)-
\log\left(\sum_{\sb\in\sper_{n+1}(X_m)}\exp(S_n\phi(\b))\right)\right)\\
               &  &\hskip 60pt \times
\exp\left(\pm \left(3\ell||\phi||\log(\#A)+5\Lambda\right)\right)
\end{eqnarray*}
which hold for any $n\leq p$, $\a\in\LL_{(m,n)}$, and $\a^*\in[\a]$. 
On the other hand, the computations performed in the proof 
of Lemma~\ref{lemma-pressure} lead us to the 
inequalities
\[
\log\left(\sum_{\sb\in\sper_{n+1}(X_m)}\exp(S_n\phi(\b))\right)=
(n+1)\left(P(\phi,X_m)\pm\left(C_\EE\theta_\EE^{n+1}+\Lambda\theta^{n+1}\right)\right),
\]
for each $m\leq n$ and $n$ such that $(n+1)(n+\ell+1)\leq p$. Since 
$\EE_{(m,p+1)}[\a]=\exp(\pm 2(p+2)C_\EE\theta^{\sqrt{p}-(\ell/2+1)})\EE_{(m,p)}[\a]$,
it follows by induction that 
$$
\mu^m[\a]=\EE_{(m,p)}[\a]\times\exp\left(
\pm 2C_\EE\theta_\EE^{-(\ell/2+1)}Q(\sqrt{p})\theta_\EE^{\sqrt{p}}\right)
$$
for each $\a\in \cup_{k=1}^{\lfloor \sqrt{p}-(\ell/2+1)\rfloor}\LL_{(m,k)}$. Therefore,
for each $m\leq n$, $\a\in\LL_{(m,n)}$, and 
$\a^*\in[\a]$, we have
\[
\frac{\mu^m[\a]}{\exp\left(S_n\phi(\a^*)-(n+1)P(\phi,X_m)\right)}=\exp(\pm C_{\rm FT}),
\]
with 
\begin{eqnarray*}
C_{\rm FT}&:=&2C_\EE\theta_\EE^{-(\ell/2+1)}
\max\left\{Q(k)\theta_\EE^k:\ k\in\nn\right\}+3\ell\parallel\!\phi\!\parallel\log(\#A)\\
                           &   &
\hskip 20pt +5\Lambda+\max\left\{(n+1)
\left(C_\EE\theta_\EE^{n+1}+\Lambda\theta^{n+1}\right):\
n\in\nn\right\}.
\end{eqnarray*}
Now, for $\a\in \LL_n$ with $n\leq m$, we obtain
$$
\mu^m[\a] = \sum_{\sb\in\LL_{m,n+k}\cap[\sa]}\mu^{m}[\b]=
$$
$$         
\exp\left(\pm C_{\rm FT}\right)\sum_{\sb\in\LL_{m,n+k}\cap[\sa]}
\exp\left(S_{n+k}\phi(\b^*)-(n+k+1)P(\phi,X_m)\right)=
$$
$$
\exp\left[S_n\phi(\a*)-(n+1)P(\phi,X_m)
\pm (C_{\rm FT}+\ell(\log(\#A)+\parallel\!\phi\!\parallel)+\Lambda)\right]
$$
$$
\times
$$
$$
\sum_{\sb\in\LL_{m,k-1}}e^{S_{k-1}\phi(\sb^*)-kP(\phi,X_m)}=
$$
$$
\exp\left[S_n\phi(\a^*)-(n+1)P(\phi,X_m)
\pm (2C_{\rm FT}+\ell(\log(\#A)+\!\parallel\!\!\phi\!\!\parallel)+\Lambda)\right]
\!\!\sum_{\sb\in\LL_{m,k-1}}\!\!\!\!\!\mu^m[\b]=
$$
$$
\exp(S_n\phi(\a*)-(n+1)P(\phi,X_m))\exp(\pm C_g)
$$
by using the specification property, and for $k$ sufficiently large. Here 
\[
C_g:=2C_{\rm FT}+\ell(\log(\#A)+\!\parallel\!\phi\!\parallel)+\Lambda). 
\]

\ms In this way we prove that
$\mu^m$ satisfies the Gibbs inequality. Theorem 1.16
in~\cite{bowenbook}, establishing the existence and uniqueness
of the Gibbs measure $\mu_\phi^m\in\MM_T(X_m)$, implies that $\mu^m:=\mu_\phi^m$.

\bs Let  
$\tilde{m}=\min\left\{m\in\nn:\ 4((k+\ell+1)^2+1)C_X\theta_X^k \leq 1 \
\text{for all}\  k\geq m\right\}$,
and $m_X$ as in Lemma~\ref{lemma-elementaryelementary}.
From Lemma~\ref{lemma-atomic}, Lemma~\ref{lemma-elementaryelementary},
and following the computations in the first part of this proof, we obtain 
$$
D(\EE_{(m,(m+\ell+1)^2)},\EE_{(m+1,(m+\ell+2)^2)}) \leq
$$
$$ 
D(\EE_{(m,(m+\ell+1)^2)},\EE_{(m,(m+\ell+2)^2)}) +
D(\EE_{(m,(m+\ell+2)^2)},\EE_{(m+1,(m+\ell+2)^2)})\leq
$$
$$
4((m+\ell+2)^2+1)C_{\EE}\theta_{\EE}^{m+1}
+2^{-m}+8((m+\ell+2)^2+1)C_X\theta_X^m,
$$
for all $m\geq m^*$, with $m^*=\max(m_X,\tilde{m})$. 

\ms Since
$\sum_{m=0}^{\infty}(m+\ell+2)^2\max(\theta_{\EE},\theta_X)^{m}$ is finite, then 
$\mu^*:=\lim_{m\to\infty}\EE_{(m,(m+\ell+1)^2)}$ is a well defined
measure in $\MM_T(X)$.
Furthermore, the convergence is such that 
\[
D(\mu^*,\EE_{(m,(m+\ell+1)^2)})\leq 2^{-m+1}+
4C_\EE Q_\EE(m)\theta_\EE^{m}+8C_X Q_X(m)\theta_X^{m-1}
\]
with 
\begin{eqnarray*}
Q_\EE(x)&:=& -\frac{(x+\ell+2)^2+1}{\log(\theta_\EE)}+
                                    \frac{2(x+\ell+2)}{\log^2(\theta_\EE)}-
                                         \frac{2}{\log^3(\theta_\EE)}\\
Q_X(x)&:=&
-\frac{(x+\ell+2)^2+1}{\log(\theta_X)}+\frac{2(x+\ell+2)}{\log^2(\theta_X)}-
                                         \frac{2}{\log^3(\theta_X)}\, .
\end{eqnarray*}

\ms Therefore, for any $m\geq m^*$, one has
$$
D(\mu^*,\mu_\phi^m)\leq D(\mu^*,\EE_{(m,(m+\ell+1)^2)})+
D(\mu_\phi^m,\EE_{(m,(m+\ell+1)^2)})\leq  Q_{\rm FT}(m)\theta_{\rm FT}^m
$$
with
$$
Q_{\rm FT}(m):=
$$
$$
4C_\EE (\theta_\EE^{-(\ell/2+1)}Q(m)+\theta_\EE^{-1}Q_\EE(m))
+8C_X\theta_X^{-1}Q_X(m)+(m+3)2^{(\ell/2+3)}+2
$$
and
$$
\theta_{\rm FT}:=\max(\theta_\EE,\theta_X,1/2)\,.
$$
\endpf

\ms
\begin{remark} In the previous proof, the polynomials
$Q,Q_\EE$, and $Q_X$
were obtained by upper bounding the series
$\sum_{k=m}^\infty P(k)\eta^k$, with $P(x)$ an increasing
polynomial, and $\eta\in(0,1)$, by the integral 
$\eta^{-1}\int_m^\infty P(x)\eta^x\ dx$.
Then we used the identity
$$
\int_m^\infty P(x)\eta^x\ dx=\eta^m\times\sum_{k=0}^{{\rm deg}(P)}
\left(-1/\log(\eta)\right)^{k+1}P^{(k)}(m)
$$
where $P^{(k)}$ is the $k$th derivative of $P$. 

\end{remark}

\bs \subsection{Proof of Theorem~\ref{theorem-gibbs-inequality}}

In the previous proof we derived the inequalities
\[\frac{\mu^m_\phi[\a]}{\exp\left(S_n\phi(\a^*)-(n+1)P(\phi,X_m)\right)}=
\exp(\pm C_g),\]
valid for each $n\in\nn$, $\a\in\LL_{(m,n)}$ and $\a^*\in[\a]$. 

On the other hand, Lemma~\ref{lemma-pressure} ensures that
$P(\phi,X_m)=P(\phi,X)\pm C_P\theta_P^m$, therefore 
\[
\frac{\mu^m_\phi[\a]}{\exp\left(S_n\phi(\a^*)-(n+1)P(\phi,X_m)\right)}=
\exp\left(\pm \left(C_g+(n+1)C_P\theta_P^m\right)\right),
\]
valid for each $n\in\nn$, $\a\in\LL_{(m,n)}$ and $\a^*\in[\a]$. Taking the limit 
$m\to\infty$, we obtain the desired result.
\endpf

\bs 
\subsection{Proof of Theorem~\ref{theorem-mixing}}

Proceeding as in the proof of Lemma~\ref{lemma-elementaryelementary}, 
the specification property implies 
$$
\frac{\di \EE_{(m,(m+\ell+1)^2)}[\a]}{\di
  \EE_{(m+1,(m+\ell+2)^2)}[\a]}=
$$
$$
\frac{\di \sum_{\sb\in\sper_{(m+\ell+2)^2+1}(X_{m+1})}\exp(S_p\phi(\b))
    }{\di \sum_{\sb\in\sper_{(m+\ell+1)^2+1}(X_{m})}\exp(S_p\phi(\b))}
\times
\frac{\di \sum_{\sb\in\sper_{(m+\ell+1)^2+1}(X_{m})\cap[\sa]}\exp(S_p\phi(\b))}{
     \di
     \sum_{\sb\in\sper_{(m+\ell+2)^2+1}(X_{m+1})\cap[\sa]}\exp(S_p\phi(\b))}=
$$
$$
\exp\left(\pm\left(4((m+\ell+2)^2+1)C_X\theta_X^m+
4\ell(\log(\#A)+||\phi||)+4\Lambda\right)\right),
$$
for each $n\leq m\in\nn$ and $\a\in\LL_{(m,n)}\equiv \LL_n$, as long as $m\geq m^*$.
These inequalities can be viewed as extensions to cylinders of 
the inequalities of Lemma~\ref{lemma-elementaryelementary}. 

\ms On the other hand, Lemma~\ref{lemma-elementarymatrix} ensures that
\[
\EE_{(m,(m+\ell+1)^2)}[\a]=\EE_{(m,(m+\ell+2)^2)}[\a]\times
\exp(\pm 2((m+\ell+2)^2+1)C_\EE\theta_\EE^{m+1}).
\] 
These and the previous inequalities imply that 
$\mu^*[\a] = \EE_{(m,(m+\ell+1)^2)}[\a]\exp(\pm\gamma_{\rm FT})$ for each $m\geq m^*$,
$m\geq n$, and $\a\in\LL_{n}$. 
Here 
\begin{eqnarray*}
\gamma_{\rm FT}&:=&4\ell(\log(\#A)+||\phi||)+4\Lambda
+\sum_{k=m^*}^\infty ((k+1)^3+1)(4C_X\theta_X^k+2C_\EE\theta_\EE^{k+1})\\
        & =&4\ell(\log(\#A)+||\phi||)+4\Lambda
+4C_XQ_X(m^*)\ \theta_X^{m^*-1}+2C_\EE Q_\EE(m^*)\theta^{m^*}.
\end{eqnarray*}

\ms Because of the previous inequalities, 
\[
\left|\mu^*\left([\a]\cap T^{-s}[\b]\right)-\mu^*[\a]\mu^*[\b]\right|
\leq 
\]
\[
e^{2\gamma_{\rm FT}}\left|\EE_{(m,(m+\ell+1)^2)}
\left([\a]\cap T^{-s}[\b]\right)-
\EE_{(m,(m+\ell+1)^2)}[\a]\EE_{(m,(m+\ell+1)^2)}[\b]\right|,
\]
for every $\a\in\LL_n$, $\b\in\LL_{n'}$ and $s\in \nn$, as long as
$n+n'+s\leq m$.

\ms Fix $\a,\b\in\LL_{m}$, $p=p(m):=(m+\ell+1)^2$, and $s\leq s'$ such that
$s+s'+4(m+1)=p+1$. Following the computations of the proof of 
Lemma~\ref{lemma-elementarymatrix}, we obtain
\begin{eqnarray*}
\EE_{(m,p)}\left([\a]\cap T^{-s}[\b]\right)&=&
\frac{\di \sum_{\sc\in\sper_{p+1}(X_m):\ \sc(0:n)=\sa,\
    \sc(n+s:n+s+n')=\sb}
\exp(S_{p}\phi(\c))}{
      \di \sum_{\sc\in\sper_{p+1}(X_m)}\exp(S_{p+1}\phi(\c))} \\
                                           &=&
\frac{\di L_{\sa}^\dag M_{(m,m)}^{s-2\ell}R^{\sb}
\times L_{\sb}^\dag M_{(m,m)}^{s'-2\ell}R^{\sa}}{
\sum_{\sc\in\LL_{(m,m)}}L_{\sc}^\dag M^{p-2\ell}R^{\sa}}
\times\exp\left(\pm 3(p+1)C\theta^{k+1}\right)
\end{eqnarray*}
$$
=\v_{(m,m)}(\a)\w_{(m,m)}(\b)\times\w_{(m,m)}(\a)\v_{(m,m)}(\b)
\times\left(\pm\left(3(p+1)C_\EE\theta^{\lfloor s/(m+\ell+1)\rfloor}\right)\right).
$$
Therefore, by using Lemma~\ref{lemma-elementarymatrix} we obtain
$$
\EE_{(m,p)}\left([\a]\cap T^{-s}[\b]\right) =
\EE_{(m,p)}[\a]\EE_{(m,p)}[\b]\times
\left(\pm\left(5(p+1)C_\EE\theta_\EE^{\lfloor s/(m+\ell+1)\rfloor}\right)\right)=
$$
$$
\EE_{(m,p)}[\a]\EE_{(m,p)}[\b]\times
\left(\pm\left(5((m+\ell+1)^2+1)
C_\EE\theta_\EE^{\lfloor s/(m+\ell+1)\rfloor}\right)\right),
$$
for each $\a,\b\in\LL_m$. Because of the additivity of the measure $\EE_{(m,p)}$, 
these inequalities extend to any $\a,\b\in\cup_{k=0}^m\LL_k$.

\ms Finally, combining the previous inequalities we obtain
\[
\left|\mu^*\left([\a]\cap T^{-s}[\b]\right)-\mu^*[\a]\mu^*[\b]\right|
\leq 
\]
\[e^{4\gamma_{\rm FT}}\left(
\exp\left(5((m+\ell+1)^2+1)C_\EE\theta_\EE^{\lfloor s/(m+\ell+1)\rfloor}\right)-1\right)
\times \mu^*[\a]\mu^*[\b],
\]
for each $\a,\b\in\cup_{k=0}^m\LL_k$. 
Let
$$
s_0:=\min\{s\in\nn:\ 
5((2\sqrt{k}+\ell+1)^2+1)C_\EE\theta_\EE^{\sqrt{k}/2-(\ell+5)/4}\leq 1\ \text{for all}
\ k\geq s\}\, .
$$
The result follows by taking 
$m=m(s):=\lfloor2\sqrt{s}\rfloor$, so that 
$$
\left|\mu^*\left([\a]\cap T^{-s}[\b]\right)-\mu^*[\a]\mu^*[\b]\right|
\leq 
$$
$$
e^{4\gamma_{\rm FT}}
\left(
\exp\left(5((2\sqrt{s}+\ell+1)^2+1)C_\EE\theta_\EE^{\sqrt{s}/2-(\ell+5)/4}\right)
-1\right) 
\times \mu^*[\a]\mu^*[\b]\leq
$$
$$
10\ C_\EE e^{4\gamma_{\rm FT}}\theta_\EE^{-(\ell+5)/4}((2\sqrt{s}+\ell+1)^2+1)
\ \theta_\EE^{\sqrt{s}/2} \times \mu^*[\a]\mu^*[\b]
$$
for all $\a\in\LL_n$, $\b\in\LL_{n'}$, and 
$s > (\max(n,n')+(\ell+1))^2/4$.
The theorem follows with
\begin{eqnarray*}
\theta_\mu&:=&\sqrt{\theta_\EE}\\
s^*(\a,\b)&:=&\max(\max(n,n')^2/4,s_0)\\
Q_\mu(x)&:=&10\ C_\EE e^{4\gamma_{\rm FT}}\ \theta_\EE^{-(\ell+5)/4}((2x+\ell+1)^2+1).
\end{eqnarray*}

\endpf

\bs \subsection{Proof of Theorem~\ref{entropy-approximation}}

By \cite{bowenbook}, each measure $\mu_{\phi}^m$ satisfies the
variational principle, as well as the measure $\mu_\phi$ by
\cite{bowen75}. This means in particular the following:
\begin{equation}\label{VP}
P(\phi,X_m)= \int_{X_m} \phi\ d\mu_\phi^m + h(\mu_\phi^m)
\quad\textup{and}\quad
P(\phi,X)= \int_{X} \phi\ d\mu_\phi + h(\mu_\phi)\,.
\end{equation}

Hence we have
$$
|h(\mu_\phi^m) - h(\mu_\phi) | \leq |P(\phi,X_m)-P(\phi,X)| + 
\left| \int_{X}\phi\ d\mu_\phi -  \int_{X_m}\phi\ d\mu_\phi^m \right|\,.
$$
It is obvious from Lemma \ref{lemma-pressure} that
\begin{equation}\label{deltaP}
0<P(\phi,X_m)-P(\phi,X) \leq \frac{C_P}{1-\theta_P} \theta_P^m\,.
\end{equation}
On another hand,
$$
\left| \int_{X}\phi\ d\mu_\phi -  \int_{X_m}\phi\ d\mu_\phi^m
\right|\leq C \theta^m\,.
$$
Statement \eqref{entropy} is thus proved.

Now, applying \eqref{relative-entropy-formula} (see appendix below)
and using \eqref{VP}-\eqref{deltaP} we get:
$$
h(\mu_\phi | \mu_\phi^m)= P(\phi,X_m)-\left(\int_X \phi\ d\mu_\phi +
h(\mu_\phi)\right)=P(\phi,X_m)-P(\phi,X)\leq \frac{C_P}{1-\theta_P} \theta_P^m\,.
$$
This proves \eqref{relative-entropy}. The proof of the theorem is now complete.

\endpf

\bs\section{Examples, Generalizations and Comments}\label{final-remarks} 

\ms A natural class of specified sofic subshifts is provided by
$\beta$-shifts coding the dynamics of the map on the unit interval
$T_\beta: x\mapsto \beta x$ mod $1$, where $\beta>1$ is a real number.
For certain $\beta$'s, the corresponding $\beta$-shift is a specified
sofic subshift. In \cite{leborgne}, the authors constructs a 
sofic coding of hyperbolic automorphisms of the torus.
In both cases, the Lebesgue measure on the unit interval or the torus is sent to the measure of maximal
entropy on the coding subshift.

\ms In this paper we assumed, for the sake of definiteness, that the potential $\phi$ was H{\"o}lder continuous 
and the subshift $X\subset A^\nn$ was a specified sofic subshift. Nevertheless, both assumptions can be weakened. 
In the proof of Theorem~\ref{theorem-elementary-convergence}, and in all other computations, the exponential decay 
\[\max\{|\phi(\a)-\phi(\b)|: \a(0:m)=\b(0:m) \}|\leq C \theta^m\] can be replaced by a polynomial decay 
\[\max\{|\phi(\a)-\phi(\b)|:  \a(0:m)=\b(0:m) \}|\leq C m^{-\overline{\alpha}},\] as long as $\overline{\alpha}>4$.
By doing so, the speed of convergence of topological pressure (Lemma
\ref{lemma-pressure}) become polynomial as well. Hence, the speed of convergence in Theorem
\ref{entropy-approximation} also become polynomial (see the proof).

\ms Regarding the nature of the subshift, the reader can verify that the essential assumptions are specification 
and presence of magic words. Moreover, the latter assumption is only
used in Lemma \ref{lemma-elementaryelementary}.
Specified sofic subshifts form a natural class of subshifts having the specification property as well as magic words, but
there are huge classes of non--sofic specified subshifts with magic words. Among them, we can mention the class of
non--sofic specified $\beta$--shifts (see~\cite{jorg}). One can
straightforwardly prove that for each non-sofic 
specified  $\beta$--shifts there exists $k\in \nn$ such that $0^k$ is a magic word. 

\ms On the other hand,
following the examples in~\cite{denker90} we can obtain non--sofic
specified subshifts with magic words, as finitary codings of
Bernoulli shifts. Take for example the finitary coding $\pi:\{0,1,2,3\}^\nn\to\{0,1,2,3\}^\nn$ such that
\[(\pi\a)_n=\left\{\begin{array}{ll} 
0     & \text{ if } \a(0:2k+1)=32^k1^k0 \text{ for some } k\in\nn\\
\a(n) & \text{ otherwise}\end{array}\right.\]  
The image subshift $X:=\overline{\pi \{0,1,2,3\}^\nn}$ is not sofic:
its description involves a non--regular language. 
Nevertheless it has the specification property, we may connect any two
admissible words by words of the kind $12^\ell1$,
and 3 is magic letter. 
Any product measure on $\{0,1,2,3\}^\nn$ induces a Gibbs measure in $X$, which can be approximated by our method. 

\ms Though the class of systems considered here is only a subclass of
those covered by Theorem 2.5 in \cite{hr92},
we are able to obtain a speed of convergence (in the weak distance) of
finite type approximations to the Gibbs measure
on the approximated subshift $X$. We were also able to prove a strong
mixing property, implying Bernoullitcity. Finally, we provide a speed of
convergence of the entropy of the finite-type approximations to the
entropy of the Gibbs measure on $X$.
We also emphasize that all constants appearing in the statements of Section \ref{main-results}
have explicit expressions in terms of the data of the problem. We did
not write these explicit formulas in the statements
because they are cumbersome. They of course appear in the course of
the proofs. 
It is also worth to notice that we only used classical algebraic tools
and symbolic dynamics, except for uniqueness of $\mu_\phi$ for which
we used Bowen's argument.

\ms  Further work has to be done in order to generalize our results to
more general subshifts. One possible approach
requires the a precise control of the convergence of the pressure. A similar approach was already exploited by Gurevich  
in the proof of the uniqueness of the maximal measure for a class of
non--specified subshift~\cite{gurevich}.
Unfortunately, the systems satisfying the hypotheses of Gurevich's theorem cannot be explicitly characterized. 

\bs{\bf Acknowledgment}. We thank K. Petersen for providing us reference \cite{gurevich}.

\bs\section{Appendix}

\bs\subsection{Primitive matrices}
$M:\{1,2,\ldots,n\}\times\{1,2,\ldots,n\}\to[0,\infty)$ is said to be
{\it primitive} if there exists an integer $\ell\geq 1$ such that
$M^\ell>0$. The smallest such integer is the {\it primitivity index}
of $M$.

\ms For $M$ primitive let
\begin{equation}\label{defintion-gamma} 
\Gamma(M):=\left\{
\begin{array}{cr}
    {\di \min_{i,j,k,l}\sqrt{\frac{M(i,j)M(k,l)}{M(i,l)M(k,j)}}} &
    M>0,\\ 
0 & \text{ otherwise.}
\end{array}\right. 
\end{equation} 
The {\it Birkhoff's coefficient} for $M$ is
$\tau(M):=(1-\Gamma(M))/(1+\Gamma(M))$.

\bs Consider the function $d:(\rr^+)^n
\times\left(\rr^+\right)^n\to\rr^+$ such that
\begin{equation}\label{definition-projectivedistance}
d(\x,\y)=\log\left(\frac{\max_i\x(i)/\y(j)}{\min_i\x(i)/\y(i)}\right).
\end{equation} 
It is the projective distance when restricted to the simplex
\[
\Delta_n:=\left\{\x:\{1,2,\ldots,n\}\to (0,1):\ |\x|_1:=
\sum_{i=1}^n\x(i)=1 \right\}.
\]

\ms The Birkhoff's coefficient gives the contraction rate of the
action of $M$ over the vector in $\Delta_n$.

\begin{tmma}\label{theorem-seneta} 
With $M$, $\Delta_n$ and $d$ be as above, define
$F_M:\Delta_n\to\Delta_n$ be such that
\[
F_M\x:=\frac{M\x}{|M\x|_1}.
\] 
Then $F_M$ is a contraction in $(\Delta_n,d)$ with contraction
coefficient $\tau(M)$, i.~e.,
\[
d(F_M\x,F_M\y)\leq \tau(M)d(\x,\y),\ \forall\ \x,\y\in\Delta_n.
\] 
\end{tmma}

\ms A proof of this result can be easily derived from the {Theorem
  3.12} in~\cite[p. 108]{senetabook}.

\ms The previous result directly implies the Perron-Frobenius Theorem
(see~\cite[ch.  1]{senetabook} for more details): a primitive matrix
$M$ has only one maximal eigenvalue $\rho_M>0$. Associated to it there
is a unique right eigenvector $\v_M\in\Delta_n$, and a unique left
eigenvector $\w_M >0$ such that $\w_M^\dag\v_M=1$.

\ms A rather direct consequence of the previous theorem is the following.

\begin{corollary}\label{corollary-projectivedistance} 
For $M$ primitive with primitivity index $\ell$, let $F:=F_{M}$ and
$\tau:=\tau(M^\ell)$. Then, for each $\x\in\Delta_n$ and $m\in\nn$, we
have
\[ 
d(F^{m}\x,\v_M)\leq \frac{\tau^{\lfloor m/\ell\rfloor}}{1-\tau}\times
d_M(\x)
\] 
with $d_M(\x):= \min\left(\ell\ d(\x,F\x),d(\x,F^{\ell}\x)\right)$.
\end{corollary}

\ms From this we readily deduce the following.

\begin{corollary}\label{corollary-powermatrix} 
Let $M:\{1,2,\ldots,n\}\times\{1,2,\ldots,n\}\to \rr^+$ be a primitive
matrix with primitivity index $\ell$, $F:=F_M$, and
$\tau:=\tau(M^\ell)$. Then, for each $\x\in\Delta_n$ and $m\in\nn$ we
have
\[
M^{m}\x=\rho_M^m\left(\w_M^{\dag}\x\right)\v_M \exp\left(
\pm\frac{\tau^{\lfloor m/\ell\rfloor} d_M(\x)}{1-\tau} \right)
\] 
with $d_M(\x):= \min\left(\ell\ d(\x,F\x),d(\x,F^{\ell}\x)\right)$.
\end{corollary}

\ms
\begin{proof} Since $M^m\x=|M^m\x|_1F^m\x$, then 
\[ 
d(F^m\x,\v_M)=\log\left(\frac{\max_i(M^m\x)(i)/\v_M(i)}{\min_i(M^m\x)(i)/\v_M(i)}\right).
\] 

\ms With
\[ 
C_m(\x):=\left(\max_i\frac{(M^m\x)(i)}{\v_M(i)}\min_i\frac{(M^m\x)(i)}{\v_M(i)}\right)^{1/2}
\] 
we have $M^m\x=C_m(\x)\v_M\times e^{\pm
d(F^m\x,\v_M)/2}$. Multiplying from the left these inequalities by
$\w_M^{\dag}$ yields $C_m(\x)=\rho_M^{m}(\w_M^{\dag}\x)e^{\pm
d(F^m\x,\v_M)/2}$.  Taking into account
Corollary~\ref{corollary-projectivedistance}, the desired result
follows.
\end{proof}

\bs\subsection{Weak distance} In this subsection $X\subset A^{\nn}$ is
any subshift. We have the following lemmas.

\begin{lemma}\label{lemma-marginalsratio} Let $\nu, \mu\in \MM(X)$ be
such that $\mu[\a]=\nu[\a]\ \exp(\pm \epsilon)$ for each
$\a\in\LL_k(X)$, then $D(\mu,\nu)\leq (\exp(\epsilon)-1)+2^{-k}$.
\end{lemma}

\begin{proof} For $j\leq k$ we have \begin{eqnarray*}
\sum_{\sa\in\LL_j(X)}|\mu[\a]-\nu[\a]|&\leq& \sum_{\sa\in\LL_j(X)}
\left(\sum_{\sb\in\LL_k(X):\ \sb(0:j)=\sa}|\mu[\b]-\nu[\b]|\right)\\
&\leq&\sum_{\sa\in\LL_j(X)} \left(\sum_{\sb\in\LL_k(X):\
\sb(0:j)=\sa}\mu[\b](e^{\epsilon}-1) \right)=e^{\epsilon}-1\, .
\end{eqnarray*}
Hence
$D(\mu,\nu)\leq (e^{\epsilon}-1)\sum_{j=0}^k2^{-(j+1)}+ \sum_{j=k+1}^\infty
2^{-(j+1)}\left(\sum_{\sa\in\LL_j(X)}|\mu[\a]-\nu[\a]|\right)$.  The
result follows taking into account that
$\sum_{\sa\in\LL_j(X)}|\mu[\a]-\nu[\a]|\leq 2$ for all $j\in\nn$.
\end{proof}

\bs \begin{lemma}\label{lemma-atomic} Let $\nu, \mu\in
\MM(X)$ be atomic with support $S_\nu:={\rm supp}(\nu)\subset
S_\mu:={\rm supp}(\mu)$.  Suppose that $\mu\{x\}\leq\nu\{x\}\leq
\mu\{x\}\exp(\epsilon)$ for each $x\in S_\nu$, then $D(\mu,\nu)\leq
\exp(\epsilon)-\exp(-\epsilon)$.
\end{lemma}

\begin{proof} For each $k\in\nn$, since $\{[\a]:\ \a\in\LL_k(X)\}$
is a partition of $X$, we have
\begin{eqnarray*} 
\sum_{\sa\in \LL_k(X)}\left|\mu[\a]-\nu[\a]\right| & = & \sum_{\sa\in
                      \LL_k(X)}(\nu(S_\nu\cap[\a])-\mu(S_\nu\cap[\a]))
                      + \sum_{\sa\in A_m}\mu((S_\mu\setminus
                      S_\nu)\cap[\a])\\ & \leq &
                      (e^{\epsilon}-1)\mu(S_\nu)+\mu(S_\mu\setminus
                      S_\nu) \leq (e^{\epsilon}-1)+\mu(S_\mu \setminus
                      S_\nu)\, .
\end{eqnarray*} 
Now, $1=\nu(S_\nu)\leq \exp(\epsilon)\mu(S_\nu)$, hence
$\mu(S_\mu\setminus S_\nu)\leq 1-\exp(-\epsilon)$ and the result
follows.
\end{proof}

\bs\subsection{Entropy and relative entropy}

Let $\nu$ be a shift-invariant probability measure on a specified subshift
$Y\subset A^{\nn}$. The measure-theoretic entropy of $\nu$ is
$$
h(\nu)=-\lim_{n\to\infty}\frac{1}{n+1} \sum_{\a\in \LL_n(Y)} \nu[\a]\log\nu[\a]\,.
$$
Since $\nu(A^{\nn}\backslash Y)=0$, we can replace $\LL_n(Y)$ by $A^{n+1}$ by
using the usual convention `$0\log 0=0'$.

We now turn to relative entropy. 
We refer the reader to \cite{CFL} for details. Therein, only subshifts
of finite type are considered but the extension to more general
subshifts is straightforward.
Let $\mu_\psi$ be a Gibbs measure
(with H{\"o}lder continuous potential $\psi$ defined on $A^{\nn}$) on a specified subshift $Y'\supset Y$.
The relative entropy of $\nu$ with respect to $\mu_\psi$ is defined
as:
$$
h(\nu|\mu_\psi )=\lim_{n\to\infty}\frac{1}{n+1}
\sum_{\a\in \LL_n(Y')} \nu[\a] \log\frac{\nu[\a]}{\mu_\psi[\a]}\,.
$$
Notice that the hypothesis $Y\subset Y'$ is crucial to make $\nu[\a]
\log\frac{\nu[\a]}{\mu_\psi[\a]}$ well-defined.
One can prove that
\begin{equation}\label{relative-entropy-formula}
h(\nu|\mu_\psi)= P(\psi,Y')-\int_{Y} \psi \ d\nu - h(\nu)\,.
\end{equation}

We notice that this result is true whenever $\mu_\psi$ satisfies the
`Gibbs inequality' \eqref{inequality-gibbs}, $\psi$ not being necessarily H{\"o}lder continuous.


\begin{thebibliography}{99} 

\bibitem{bowen75} R. Bowen, {\em Some systems with unique equilibrium
states}, Math.  Systems Theory {\bf 8} (1974/75), no. 3, 193--202.
  
\bibitem{bowenbook} R.~Bowen, {\em Equilibrium States and the Ergodic
Theory of Anosov Diffeomorphisms}, Lecture Notes in Mathematics {\bf
470}, Springer--Verlag, 1975.

\bibitem{CFL} 
J.-R. Chazottes, E. Floriani, R. Lima,
{\em Relative entropy and identification of Gibbs measures in dynamical systems},
J. Statist. Phys. {\bf 90} (1998), no. 3-4, 697--725.

\bibitem{denkerbook} M. Denker, C. Grillenberger, K. Sigmund, {\em
Ergodic Theory on Compact Spaces}, Lecture Notes in Math. {\bf 527},
Springer-Verlag (1976).

\bibitem{denker90} M. Denker, {\em Some New Examples of Gibbs 
Measures}, Monat. fur Math. {\bf 109} (1990) 49--62.

\bibitem{gurevich} B. Gurevich, 
{\em Stationary random sequences of maximal entropy}.  In 
Multicomponent random systems,  pp. 327--380,
Adv. Probab. Related Topics {\bf 6} Dekker, New York, 1980.

\bibitem{hr92} N.~T.~A.~Haydn and D.~Ruelle, {Equivalence of Gibbs and
Equilibrium States for Homeomorphisms Satisfying Expansiveness and
Specification}, Commun. Math. Phys. {\bf 148} (1992), 155--167.
 
\bibitem{KH}
A. Katok, B. Hasselblatt, 
Introduction to the modern theory of dynamical systems. Encyclopedia
of Mathematics and its Applications {\bf 54}.
Cambridge University Press, Cambridge, 1995. 

\bibitem{kellerbook} G. Keller, {\em Equilibrium States in Ergodic
Theory}, London Mathematical Society Student Texts {\bf 42}, Cambridge
University Press 1998.
 
\bibitem{kitchensbook} B. Kitchens, {\em Symbolic Dynamics},
Springer-Verlag, Berlin, 1998.

\bibitem{leborgne}
S. Le Borgne, 
{\em Un codage sofique des automorphismes hyperboliques du tore},
S{\'e}minaires de Probabilit{\'e}s de Rennes (1995), 35 pp.,
Publ. Inst. Rech. Math. Rennes, 1995, Univ. Rennes I, Rennes, 1995. 


\bibitem{ruelle73} D.~Ruelle {\it Statistical Mechanics on compact
sets with $\zz^\nu$ actions satisfying expansiveness and
specification}, Trans. Amer. Math. Soc. {\bf 185} (1973), 237--251.

\bibitem{jorg}
J. Schmeling, 
{\em Symbolic dynamics for $\beta$-shifts and self-normal numbers},
Ergodic Theory Dynam. Systems {\bf 17} (1997), no. 3, 675--694.

\bibitem{senetabook} E.~Seneta, {\em Non-negative Matrices and Markov
Chains}, Springer Series in Statistics, Springer-Verlag, 1981.

\bibitem{shieldsbook} P.~Shields, {\em The Ergodic Theory of Discrete
  Sample Paths}, Graduate Studies in Mathematics {\bf 13}, American
  Mathematical
Society, 1996.
  
\bibitem{waltersbook} P.~Walters, {\em An introduction to Ergodic
Theory}, Springer Verlag, 1982.
 
\end{thebibliography}
\end{document}